\documentclass[12pt, reqno]{amsart}
\usepackage[active]{srcltx}
\usepackage{verbatim}
\usepackage{epsfig,graphicx,color,mathrsfs}
\usepackage{graphicx}
\usepackage{amsmath,amssymb,amsthm,amsfonts}
\usepackage{amssymb}
\usepackage[english]{babel}
\usepackage{appendix}
\usepackage{todonotes}
\usepackage[hyperpageref]{backref}
\usepackage[colorlinks=true,urlcolor=blue,citecolor=red,linkcolor=blue,linktocpage,pdfpagelabels,bookmarksnumbered,bookmarksopen]{hyperref}
\DeclareRobustCommand{\rev}[1]{{\color{red}#1}}

\usepackage[left=2.7cm,right=2.7cm,top=3cm,bottom=3cm]{geometry}

\newcommand{\N}{{\mathbb N}}

\newcommand{\R}{{\mathbb R}}

\numberwithin{equation}{section}
\newtheorem{theorem}{Theorem}[section]
\newtheorem{proposition}[theorem]{Proposition}
\newtheorem{lemma}[theorem]{Lemma}

\newtheorem{definition}[theorem]{Definition}
\newtheorem{remark}[theorem]{Remark}
\theoremstyle{definition}

\renewcommand{\dfrac}{\displaystyle\frac}
\newcommand{\brm}{\begin{remark}\rm}
\newcommand{\erm}{\end{remark}}
\newcommand{\brms}{\begin{remark}\rm}
\newcommand{\erms}{\end{remark}}
\newcommand{\bte}{\begin{theorem}}
\newcommand{\ete}{\end{theorem}}
\newcommand{\bpr}{\begin{proposition}}
\newcommand{\epr}{\end{proposition}}
\newcommand{\ble}{\begin{lemma}}
\newcommand{\ele}{\end{lemma}}
\newcommand{\beq}{\begin{equation}}
\newcommand{\eeq}{\end{equation}}
\newcommand{\bdm}{\begin{displaymath}}
\newcommand{\edm}{\end{displaymath}}
\numberwithin{equation}{section}

\newcommand{\bos}{\begin{remark}\rm}
\newcommand{\eos}{\end{remark}}

\newcommand{\ben}{\begin{enumerate}}
\newcommand{\een}{\end{enumerate}}
\newcommand{\n }{\nabla }

\newcommand{\be}{\begin{equation}}
\newcommand{\ee}{\end{equation}}

\title[Singular Finsler $p$-Laplace equations]{Singular Finsler $p$-Laplace equations in exterior domains}

\author{Francesco Esposito$^*$, Luigi Muglia$^*$ and Alessandro Trombetta$^*$}

\address{$^{*}$Dipartimento di Matematica e Informatica, UNICAL, Ponte Pietro  Bucci 31B, 87036 Arcavacata di Rende, Cosenza, Italy}

\email{francesco.esposito@unical.it, luigi.muglia@unical.it, alessandro.trombetta@unical.it}



\thanks{The authors are supported by the Gruppo Nazionale per l'Analisi Matematica, la Probabilità e le loro Applicazioni (GNAMPA) of the Istituto Nazionale di Alta Matematica (INdAM). F. Esposito is partially supported by \emph{INdAM-GNAMPA Esistenza, regolarità e proprietà qualitative per problemi non lineari} E5324001950001. L. Muglia is partially supported by \emph{INdAM-GNAMPA Esistenza e proprietà qualitative per soluzioni di equazioni ellittiche con termini di ordine inferiore} E5324001950001. F. Esposito is partially supported by the PID Project (Spain) PID2024-155314NB-I00.}

\begin{document}

\begin{abstract}
	The aim of this paper is to prove existence and uniqueness of weak solutions to the singular anisotropic problem
	\[
	-\Delta^H_p u=\frac{K(x)}{u^{\gamma}} \ \text{ in } \Omega^e,\qquad u>0 \ \text{ in } \Omega^e,\qquad u=0 \ \text{ on } \partial\Omega,\qquad u(x)\to0 \ \text{ as } |x|\to\infty,
	\]
	where $\Delta^H_p$ is the Finsler $p$-Laplacian, $H$ is a Finsler norm, $\Omega\subset\mathbb{R}^N$, $N\geq3$, is a smooth bounded domain, $\Omega^e:=\mathbb{R}^N\setminus\overline{\Omega}$ is the corresponding exterior domain, $1<p<N$, $\gamma>0$, and $K$ is a positive weight in $L^{\frac{N}{p^-}}(\Omega^e)\cap L^1(\Omega^e)$. We are able to cover both the strongly singular case $\gamma\geq1$ and the case $0<\gamma<1$, and we allow for unbounded weights. In this way, even in the isotropic case $H(\xi)=|\xi|$, we extend previous results. Existence is obtained by approximating the problem with non-singular variational problems. Uniqueness follows from a weak comparison principle for sub and supersolutions that belong only to $W^{1,p}_{loc}$ and need not have finite energy.
\end{abstract}

\maketitle

\section{Introduction}\label{sec:intro}

In this paper we study existence and uniqueness of weak solutions to the following singular quasilinear problem, driven by an anisotropic operator in an exterior domain:
\begin{equation}\label{problem}
	\begin{cases}
		-\Delta^H_p u=\dfrac{K(x)}{u^{\gamma}} & \text{in } \Omega^e,\\[3mm]
		u>0 & \text{in } \Omega^e,\\
		u=0 & \text{on } \partial\Omega,\\
		u(x)\to 0 & \text{as } |x|\to\infty .
	\end{cases}
\end{equation}
Here $\Omega\subset\mathbb{R}^N$, $N\geq 3$, is a simply connected bounded domain with smooth boundary and $0 \in \Omega$, $\Omega^e:=\mathbb{R}^N\setminus\overline{\Omega}$ is the corresponding exterior domain, $1<p<N$, $\gamma>0$, and $K:\Omega^e\to\mathbb{R}$ is a positive weight satisfying
\begin{equation}\label{hyp:K}
	K\in L^{\frac{N}{p^-}}(\Omega^e)\cap L^1(\Omega^e),
\end{equation}
where, here and throughout the paper, $K\in L^{\frac{N}{p^-}}(\Omega^e)$ means that $K\in L^{\frac{N}{p-\varepsilon}}(\Omega^e)$ for every $\varepsilon\in(0,p)$. Since $\frac{N}{p-\varepsilon}>\frac{N}{p}>(p^*)'$, where $p^*:=\frac{Np}{N-p}$ is the critical Sobolev exponent, by interpolation, from \eqref{hyp:K} we deduce that $K\in L^{(p^*)'}(\Omega^e)$.

The nonlinear operator
\begin{equation}\label{operator}
	\Delta^H_p u:=\operatorname{div}\left(H^{p-1}(\nabla u)\nabla H(\nabla u)\right)
\end{equation}
is the so-called Finsler $p$-Laplacian (or anisotropic $p$-Laplacian), namely, $-\Delta^H_p$ is (in the weak sense) the G\^ateaux derivative of the functional $u\mapsto\frac1p\int H^p(\nabla u)\,dx$, where $H$ is a Finsler norm satisfying the following hypotheses:
\begin{itemize}
	\item[(i)] $H\in C^2_{loc}(\mathbb{R}^N\setminus\{0\})$ is even and $H(\xi)>0$ for every $\xi\in\mathbb{R}^N\setminus\{0\}$;
	\item[(ii)] $H(t\xi)=tH(\xi)$ for every $\xi\in\mathbb{R}^N\setminus\{0\}$ and every $t>0$;
	\item[(iii)] $H$ is uniformly elliptic, that is, the set $B^H_1:=\{\xi\in\mathbb{R}^N:\ H(\xi)<1\}$ is uniformly convex:
	\[
	\exists\,\lambda>0:\quad \langle D^2H(\xi)v,v\rangle\geq\lambda|v|^2 \qquad \forall\,\xi\in\partial B^H_1,\ \forall\,v\in\nabla H(\xi)^{\perp}.
	\]
\end{itemize}
We refer to \cite{CFV} for more details on these assumptions. Since $H$ is a norm, there exist $a,b>0$ such that $a|\xi|\leq H(\xi)\leq b|\xi|$ for every $\xi\in\mathbb{R}^N$. When $H(\xi)=|\xi|$, the operator $\Delta^H_p$ reduces to the classical $p$-Laplacian $\Delta_p$, while, for instance, $H(\xi)=\langle A\xi,\xi\rangle^{1/2}$, with $A$ a symmetric positive definite matrix, gives rise to an anisotropic operator.

Problem \eqref{problem} combines several difficulties. First, the nonlinearity is singular, i.e. the right-hand side blows up as $u\to0^+$, and this happens near $\partial\Omega$ because of the Dirichlet condition. As a consequence, solutions are not expected to be smooth up to the boundary and their gradient may blow up near $\partial\Omega$; already for the Laplacian in bounded domains, solutions do not have finite energy when $\gamma$ is large (see the celebrated paper by Lazer and McKenna \cite{LM}). Moreover, in the strongly singular case $\gamma\geq1$ the energy functional formally associated with \eqref{problem} is not finite on the whole energy space, so that standard variational techniques cannot be applied directly. Second, the problem is set in an exterior domain: $\Omega^e$ is unbounded, the Sobolev embeddings are not compact, and the behaviour of the solutions at infinity has to be controlled in order to recover the decay condition in \eqref{problem}. Third, the operator is anisotropic: $\Delta^H_p$ is a quasilinear operator, degenerate if $p>2$ and singular if $1<p<2$, whose structure is dictated by the Finsler norm $H$. Hence, the arguments relying on the explicit form of the $p$-Laplacian, or on its invariance under rotations, are not available, and all the estimates have to be derived from the structural properties of $H$ only.

Because of the lack of regularity up to the boundary, the Dirichlet condition has to be understood in a generalized sense, following the approach of \cite{CD,CS,CST,CF,EST}. Denoting by $\mathcal{D}^{1,p}_0(\Omega^e)$ the space of functions $u\in L^{p^*}(\Omega^e)$ with $\nabla u\in L^p(\Omega^e)$ that vanish on $\partial\Omega$ in the sense of traces (the functional setting is recalled in Section~\ref{sec:Preliminaries}), we will use the following notion of solution.
\begin{definition}\label{def:sol}
	A solution of problem \eqref{problem} is a function $u\in W^{1,p}_{loc}(\Omega^e)$ such that
	\begin{itemize}
		\item[(i)] $u$ satisfies
		\begin{equation}\label{weak equation}
		\int_{\Omega^e}H^{p-1}(\nabla u)\langle\nabla H(\nabla u),\nabla\varphi\rangle\,dx=\int_{\Omega^e}\frac{K(x)}{u^{\gamma}}\,\varphi\,dx\qquad\forall\,\varphi\in C^1_c(\Omega^e);
		\end{equation}
		\item[(ii)] $u$ is strictly positive on each compact subset of $\Omega^e$, i.e. for every $\omega\subset\subset\Omega^e$ there exists $c_\omega>0$ such that $u\geq c_\omega$ in $\omega$;
		\item[(iii)] $u=0$ on $\partial\Omega$, in the sense that $u\geq0$ and $(u-\varepsilon)^+\in\mathcal{D}^{1,p}_0(\Omega^e)$ for every $\varepsilon>0$;
		\item[(iv)] $u(x)\to0$ as $|x|\to\infty$.
	\end{itemize}
\end{definition}
Condition (ii) guarantees that the right-hand side in (i) is well defined, while (iii) and (iv) express, respectively, the boundary condition on $\partial\Omega$ and the decay at infinity.

The study of elliptic equations with singular nonlinearities goes back to the seminal paper of Crandall, Rabinowitz and Tartar \cite{CRT}, where the existence of solutions was proved in the semilinear case. Subsequently, Lazer and McKenna \cite{LM}  showed that solutions are in general not smooth up to the boundary and that, for the Laplacian, they do not belong to $H^1_0$ when $\gamma\geq3$. In recent years, singular problems have received considerable attention; without any claim of completeness, we refer to \cite{CES,CGS,dP,ES,FMS,GOW,HSS,LS,MMS,OP1,OP2} and to \cite{ABLP,ACLMOP,Bo} for problems involving first-order terms. For a comprehensive account of the history of the problem, we refer to the survey by Oliva and Petitta \cite{OPsurvey}. A turning point for non-regular data is the paper of Boccardo and Orsina \cite{BO}, where an approximation and truncation technique provides existence and summability of solutions depending on $\gamma$ and on the summability of the datum; this approach was extended to the $p$-Laplacian in \cite{CST,DeC}. In the strongly singular regime, existence results have also been obtained via sub and supersolutions and truncation \cite{PS}, and via nonsmooth critical point theory \cite{CD,HSS2}; in particular, the generalized boundary condition in Definition~\ref{def:sol} goes back to \cite{CD}. Uniqueness is a more subtle issue; it is expected, since the nonlinearity is decreasing in $u$, but the lack of finite energy prevents the use of the standard arguments. A positive answer under mild assumptions on the datum was given in \cite{CS} for the Laplacian and in \cite{CST} for the $p$-Laplacian; see also \cite{DDO,OP2,OrP}. In the anisotropic setting, existence and uniqueness for the Finsler $p$-Laplacian with singular nonlinearities in bounded domains were established in \cite{EST}, under very general summability assumptions on the source term, while the question of whether solutions have finite energy has been addressed in the recent paper \cite{EOV}.

Much less is known in exterior domains. Singular quasilinear problems in $\Omega^e$ were studied in \cite{CDS,DS} by means of sub and supersolutions, for nonlinearities of the form $f(u)/u^{\gamma}$ with $0<\gamma<1$ which are $p$-sublinear at infinity. The strongly singular case was addressed for the first time by Chhetri and Faraci \cite{CF}, who proved existence and uniqueness of a weak solution to $-\Delta_p u=K(x)/u^{\gamma}$ in $\Omega^e$, with $\gamma\geq1$ and $K\in L^\infty(\Omega^e)\cap L^1(\Omega^e)$, provided decay estimates near $\partial\Omega$ and at infinity when $\Omega$ is a ball, and also considered more general nonlinearities.

The aim of the present paper is twofold: on the one hand, we extend the existence and uniqueness results of \cite{CF} to the anisotropic setting under more general assumptions; on the other hand, we carry over the results of \cite{EST} from bounded to exterior domains. Compared with \cite[Theorem~1.1]{CF}, we deal with the Finsler $p$-Laplacian, we consider every $\gamma>0$, thus covering both the strongly singular case $\gamma\geq1$ and the case $0<\gamma<1$, and we replace the assumption $K\in L^\infty(\Omega^e)\cap L^1(\Omega^e)$ with the weaker summability condition \eqref{hyp:K}, which allows for unbounded weights. In particular, even in the isotropic case $H(\xi)=|\xi|$, Theorems~\ref{thm:existence} and \ref{thm:uniqueness} below extend the existence and uniqueness part of \cite[Theorem~1.1]{CF}.

Our first main result concerns existence.

\begin{theorem}\label{thm:existence}
	Let $\gamma>0$ and let $K$ be a positive function satisfying \eqref{hyp:K}. Then problem \eqref{problem} admits a solution $u$ in the sense of Definition~\ref{def:sol}. Moreover:
	\begin{itemize}
		\item[(a)] if $\gamma\geq1$, then $u^{\frac{p+\gamma-1}{p}}\in\mathcal{D}^{1,p}_0(\Omega^e)$; in particular, $u\in\mathcal{D}^{1,p}_0(\Omega^e)$ if $\gamma=1$;
		\item[(b)] if $0<\gamma<1$, then $u\in\mathcal{D}^{1,p}_0(\Omega^e)$.
	\end{itemize}
\end{theorem}

The second main result is a weak comparison principle, which in turn yields uniqueness.

\begin{theorem}\label{thm:uniqueness}
	Let $\gamma>0$ and let $K$ be a positive function satisfying \eqref{hyp:K}. If $u,v\in W^{1,p}_{loc}(\Omega^e)$ are, respectively, a subsolution and a supersolution of \eqref{problem} in the sense of Definition~\ref{def:subsuper}, then $u\leq v$ a.e. in $\Omega^e$. In particular, problem \eqref{problem} admits a unique solution in the sense of Definition~\ref{def:sol}.
\end{theorem}

Let us briefly comment on the proofs. The existence part relies on an approximation scheme in the spirit of \cite{BO}: the singular term is regularized by means of $g_n(x,s):=K(x)\,(s^++\frac1n)^{-\gamma}$, and the corresponding non-singular problems are solved by minimizing the associated energy functionals in $\mathcal{D}^{1,p}_0(\Omega^e)$. The lack of compactness due to the unboundedness of $\Omega^e$ is overcome by splitting the domain into a bounded region and a neighbourhood of infinity, where the summability of $K$ provides the required smallness. We then show that the sequence $(u_n)$ of approximating solutions is nondecreasing, uniformly bounded away from zero on compact subsets of $\Omega^e$ and bounded in $W^{1,p}_{loc}(\Omega^e)$; this is achieved through uniform estimates on $u_n^{\frac{p+\gamma-1}{p}}$ in $\mathcal{D}^{1,p}_0(\Omega^e)$ when $\gamma\geq1$, and directly on $u_n$ in $\mathcal{D}^{1,p}_0(\Omega^e)$ when $0<\gamma<1$. Since $K$ is not assumed to be bounded, local $L^\infty$ bounds for the approximating solutions are obtained by a Moser-type iteration in the spirit of Serrin \cite{Se}, which only requires \eqref{hyp:K}. Finally, we prove that $u_n\leq u\leq u_n+\frac1n$ in $\Omega^e$, which yields the uniform convergence of $(u_n)$ and the decay of $u$ at infinity.

The proof of Theorem~\ref{thm:uniqueness} is the most delicate part of the paper, since sub and supersolutions are only locally in $W^{1,p}$ and are not required to have finite energy. Differently from \cite{CF}, where uniqueness is obtained by comparing an arbitrary solution with the approximating sequence, we follow the approach introduced in \cite{CS,CST} and adapted to the anisotropic framework in \cite{EST}. It is based on the minimization of a suitably truncated functional on the convex set of functions $\varphi\in\mathcal{D}^{1,p}_0(\Omega^e)$ with $0\leq\varphi\leq v$, combined with a careful choice of test functions, and it provides a comparison principle which we believe to be of independent interest.

We conclude by recalling that anisotropic operators of Finsler type naturally arise in Finsler geometry \cite{BCS,BP,CiS} and in several applied models, for instance in material science \cite{HC}, relativity \cite{As,CJM}, biology \cite{AIM} and image processing \cite{PM}; see also \cite{CRS, CFV} for regularity and qualitative properties of solutions to anisotropic equations.

The paper is organized as follows. In Section~\ref{sec:Preliminaries} we collect some preliminary properties of Finsler norms and we recall the functional setting. In Section~\ref{sec:existence} we prove Theorem~\ref{thm:existence}, treating separately the cases $\gamma\geq1$ and $0<\gamma<1$. Section~\ref{sec:uniqueness} is devoted to the weak comparison principle and to the proof of Theorem~\ref{thm:uniqueness}.

\section{Preliminaries}\label{sec:Preliminaries}

In this section we collect some properties of the Finsler norm $H$, which satisfies assumptions (i)--(iii) of Section~\ref{sec:intro}, and we recall the functional setting.
\begin{remark} Next we summarize some useful properties of the function $H$ proved in \cite{CFV14,CFV}.
\begin{itemize}
\item $H$ is a norm equivalent to the Euclidean one, i.e. there exist $a,b>0$ such that:
\begin{equation}\label{norm}
a|\xi|\leq H(\xi)\leq b|\xi|,\quad\forall\, \xi\in\R^N.
\end{equation}

\item $H$ is 1-homogeneous, hence, by the Euler Theorem it follows
\begin{equation}\label{eq:eulero}
 \langle\nabla H(\xi), \xi\rangle = H(\xi), \qquad\forall\, \xi \in \R^N\setminus \{0\}.
\end{equation}%\begin{align*}
Moreover
\begin{equation}\label{eq:grad 0 omog}
\nabla H(t\xi)=\hbox{sign}(t)\n H(\xi), \qquad\forall\,\xi \in \R^N\setminus \{0\},\, \forall t\neq 0.
\end{equation}
By the previous equality, we infer that there exists $C_H>0$ such that 
\begin{equation}\label{grad H limitato}
|\nabla H(\xi)| \le  C_H, \qquad \forall \xi\in\R^N\setminus \{0\}.
\end{equation}
\item Since $H^2$ is 2-homogeneous there exists $\Lambda_2>0$ such that
$$
 \langle D^2 H^2(\xi)\eta,\eta\rangle \leq \Lambda_2 |\eta|^2
$$
for $\xi\neq 0$ and $\eta\in\mathbb{R}^N$.
\item Since $D^2H$ is a $(-1)-$homogeneous function we have
\begin{equation}\label{eq:hess -1 omog}
D^2 H(t\xi)=\frac{1}{|t|}D^2H(\xi) \qquad\forall\,\xi \in  \R^N\setminus \{0\},\, \forall t\neq 0.
\end{equation}
It immediately implies that there exists $K_2 >0$ such that 
\begin{equation} \label{eq:Hess H limitato}
    |D^{2}H(\xi)| \le \frac{K_2}{|\xi|}, \qquad \forall \xi\in\R^N\setminus \{0\},
\end{equation}
where $|\cdot|$ denotes the usual Euclidean  norm of a matrix, and 
\begin{equation}\label{eq:hess nullo}
D^2 H(\xi)\xi=0, \qquad\forall\,\xi\in\R^N\setminus \{0\}.
\end{equation}

\item Assumption (iii) is equivalent to the following:
\begin{equation} \label{eq:definita positiva 2}
    \exists \Lambda > 0: \langle D^2 H(\xi)v, v \rangle \geq \Lambda |\xi|^{-1} |v|^2 \quad \forall \xi \in \R^{N} \setminus \{0\}, \; \forall v \in \nabla H(\xi)^\bot.
\end{equation}
\end{itemize}
\end{remark} \ \\
The proof of the following technical result can be found in \cite[Proposition 2.1]{EMSV} and \cite[Proposition 3.1]{CFV} taking into account $(iii)$.
\begin{lemma}\label{estimates}
	There exist $C_{1},C_{2}>0$ such that for all $\eta,\eta^{\prime}\in\mathbb{R}^N$ with $|\eta|+|\eta^{\prime}|>0$
	\begin{eqnarray}
		&&\label{1.4}\left\langle H ^{p-1} \left( \eta \right) \nabla H \left( \eta \right)-H ^{p-1} \left( \eta^{\prime}\right)\nabla H \left(\eta^{\prime} \right), \eta - \eta^{\prime} \right\rangle \geq C_{1}\left( \left\vert \eta \right\vert +\left\vert \eta^{\prime}\right\vert \right) ^{p-2}\left\vert \eta-\eta^{\prime} \right\vert ^{2}, \\ \label{1.5}&&\left\vert H ^{p-1} \left( \eta \right) \nabla H \left( \eta \right)-H ^{p-1} \left( \eta^{\prime}\right)\nabla H \left(\eta^{\prime} \right) \right\vert \leq C_{2}\left( \left\vert  \eta \right\vert +\left\vert \eta^{\prime}\right\vert \right) ^{p-2}\left\vert \eta-\eta^{\prime} \right\vert.
	\end{eqnarray}
\end{lemma}\ \\
Recall that the space
\[
\mathcal{D}_0^{1,p}(\Omega^e) :=\overline{C_c^\infty(\Omega^e)}^{\|\cdot\|},
\]
where 
\[
\|u\| := \left( \int_{\Omega^e} |\nabla u|^p \right)^{\frac{1}{p}},
\]
As usual, we denote by $L^q(\Omega^e)$ the Lebesgue space of $q$-summable functions ($1 \le q < \infty$), with the classical norm 
\[
\|u\|_q := \left( \int_{\Omega^e} |u|^q \right)^{\frac{1}{q}}.
\]
Since $1<p<N$, $\mathcal{D}_0^{1,p}(\Omega^e)$ is a reflexive Banach space which coincides with the space of functions $u\in L^{p^*}(\Omega^e)$ with $\nabla u\in L^p(\Omega^e)$ such that $u=0$ on $\partial\Omega$ in the sense of traces.\\
It is well known that $\mathcal{D}_0^{1,p}(\Omega^e)$ embeds continuously in $L^{p^*}(\Omega^e)$, where $p^* := \frac{pN}{N-p}$ is the critical Sobolev exponent, with embedding constant $C_S$.

\section{Existence of solutions to problem \eqref{problem}}\label{sec:existence}

Let $r > 0$ be such that $\overline{\Omega} \subset B_r := \{x \in \mathbb{R}^N : |x| < r\}$. We set $\Omega_r := \Omega^e \cap B_r$ and $B_r^e := \Omega^e \setminus \overline{B}_r = \mathbb{R}^N \setminus \overline{B}_r$.\ \\

\noindent Let $g_n : \Omega^e \times \mathbb{R} \to \mathbb{R}$ be defined by
\begin{equation*}
	g_n(x, s) := \frac{K(x)}{(s^+ + \frac{1}{n})^\gamma} \quad \text{for every } n \in \mathbb{N}^+,
\end{equation*}
and let $G_n : \Omega^e \times \mathbb{R} \to \mathbb{R}$ be the primitive of $g_n$ given by
\begin{equation*}
	G_n(x, s) := \int_0^s g_n(x, z) \, dz.
\end{equation*}
\begin{remark}
	It is easy to verify that $g_n$ and $G_n$ satisfy the following inequalities:
	\begin{equation}\label{stima g_n}
		0< g_n(x, s) \le n^\gamma K(x), \text{ a.e. } x \in \Omega^e, s \in \mathbb{R},
	\end{equation}
	\begin{equation}\label{stima G_n}
		|G_n(x, s)| \le n^\gamma K(x)|s|, \text{ a.e. } x \in \Omega^e, s \in \mathbb{R}.
	\end{equation}
\end{remark}

\noindent Let us consider, for every fixed $n \in \N$, the truncated problem of \eqref{problem}:
\begin{equation} \label{appr problem}
	\left\{ 
	\begin{array}{cc}
		-\Delta^{H}_{p}u=g_n(x, u) & \text{in }\Omega^e , \\ 
		u=0 & \text{on }\partial \Omega .%
	\end{array}%
	\right.
\end{equation}

\subsection{The case $\gamma\geq 1$} The next lemma proves the existence of a solution to \eqref{appr problem}. We point out that Lemma~\ref{existence appr problem}, items \emph{(i)}--\emph{(ii)} of Lemma~\ref{properties of u_n} and Lemma~\ref{lem:sandwich} below hold for every $\gamma>0$, since their proofs do not use the assumption $\gamma\geq1$.

\begin{lemma}\label{existence appr problem}
	Let $n\in \mathbb{N}$ be fixed. Problem \eqref{appr problem} admits a weak solution 
	$u_{n}\in \mathcal{D}_{0}^{1,p}\left( \Omega^e \right) $, i. e. $u_{n}$ satisfies
	\begin{equation}\label{weak equation appr problem}
		\int_{\Omega^e }H ^{p-1} \left( \nabla u_n\right) \left\langle \nabla H \left(\nabla u_n \right), \nabla \varphi \right\rangle  
		=\int_{\Omega^e } g_n(x, u_n)\varphi, \quad \forall \varphi
		\in \mathcal{D}_{0}^{1,p}\left( \Omega^e \right) .
	\end{equation}
\end{lemma}
\begin{proof}
Let $\mathcal{J}_n : \mathcal{D}_0^{1,p}(\Omega^e) \to \mathbb{R}$ defined by
\begin{equation*}
	\mathcal{J}_n(u) = \frac{1}{p} \int_{\Omega^e }H ^{p} \left( \nabla u\right) - \int_{\Omega^e} G_n(x, u).
\end{equation*}
First let us note that  $\mathcal{J}_n$ is coercive. Indeed, let $\left\lbrace w_k \right\rbrace$ in $\mathcal{D}_{0}^{1,p}\left( \Omega^e \right)$ such that $\|w_k\| \to \infty$; by \eqref{stima G_n}
\begin{equation*}
	\int_{\Omega^e} G_n(x, w_k) \leq 	\int_{\Omega^e} n^\gamma K(x)|w_k|.
\end{equation*}
Using H\"older inequality and Sobolev inequality, we obtain
\begin{equation*}
	\int_{\Omega^e} G_n(x, w_k) \leq 	\left( \int_{\Omega^e}\left(  n^\gamma K(x)\right)^{\frac{p^*}{p^*-1}}\right) ^{\frac{p^*-1}{p^*}} \left( \int_{\Omega^e} \left| w_k\right| ^{p^*}\right) ^{\frac{1}{p^*}} \leq C \|w_k\|.
\end{equation*}
where the positive constant $C=C(C_S, n, \|K\|_{L^{(p^*)'}})$.	\\ 
Hence, by \eqref{norm}
\begin{equation*}
	\mathcal{J}_n (w_k ) \geq \frac{1}{p} \int_{\Omega^e }H ^{p} \left( \nabla w_k\right) - C \|w_k\| \geq \frac{a^p}{p} \|w_k\| ^p - C \|w_k\| \to \infty,
\end{equation*}
as $\|w_k\| \to \infty$. Then $\mathcal{J}_n$ is coercive.  \\

Next let us prove that $\mathcal{J}_n$ is sequentially weakly lower semicontinuous. In order to do this, let $\left\lbrace w_k \right\rbrace$ in $\mathcal{D}_{0}^{1,p}\left( \Omega^e \right)$ such that $w_k \rightharpoonup w$; then, the sequence $\{\|w_k\|_{p^*}\}$ is bounded and $w_k \to w$ in $L_{loc}^q(\Omega^e)$ for each $1 < q < p^*$. We have
\begin{equation}
	\int_{\Omega^e} |G_n(x, w_k) - G_n(x, w)| = \int_{\Omega_r} |G_n(x, w_k) - G_n(x, w)|  + \int_{B_r^e} |G_n(x, w_k) - G_n(x, w)|.
\end{equation}
Fix $\varepsilon >0$. Observe that, for $k$ big enough, since $|G_n(x,w_k)-G_n(x,w)|\leq n^\gamma K(x)|w_k-w|$, by the H\"older inequality, the compactness of the embedding of $\mathcal{D}_0^{1,p}(\Omega^e)$ into $L^q(\Omega_r)$, $1\leq q<p^*$ and \eqref{hyp:K} we get
\begin{equation}\label{G_n 1}
	\int_{\Omega_r} |G_n(x, w_k) - G_n(x, w)| < \frac{\varepsilon}{2}.
\end{equation}
Moreover, for each $k$, by \eqref{stima g_n}, \eqref{stima G_n} and H\"older inequality, we have
\begin{align*}
	\int_{B_r^e} |G_n(x, w_k) - G_n(x, w)|
	&\le \int_{B_r^e} n^\gamma K(x)|w_k| + \int_{B_r^e} n^\gamma K(x)|w| \\
	&\le C_1 \left( \int_{B_r^e} K(x)^{(p^*)'}  \right)^{\frac{1}{(p^*)'}}(\|w_k\|_{p^*} + \|w\|_{p^*})  \\
	&\le C_2 \left( \int_{B_r^e} K(x)^{(p^*)'}  \right)^{\frac{1}{(p^*)'}},
\end{align*}
for some positive constants $C_1=C_1(n)$ and $C_2=C_2(n,\sup_k\|w_k\|_{p^*},\|w\|_{p^*})$. Choosing a suitable $r$, independent of $k$, since $K$ belongs to $L^{(p^*)'}(\Omega^e)$, we get
\begin{equation}\label{G_n 2}
	\int_{B_r^e} |G_n(x, w_k) - G_n(x, w)| < \frac{\varepsilon}{2}.
\end{equation}
Hence, choosing first $r$ as in \eqref{G_n 2} and then $k$ so large that \eqref{G_n 1} holds, by \eqref{G_n 1} and \eqref{G_n 2} we obtain
\begin{equation}
	\lim_{k \to \infty} \int_{\Omega^e} G_n(x, w_k)  = \int_{\Omega^e} G_n(x, w).
\end{equation}	

Hence, since $w\mapsto\frac1p\int_{\Omega^e}H^p(\nabla w)$ is convex and continuous, and thus weakly lower semicontinuous, $\liminf_{k \to \infty} \mathcal{J}_n (w_k ) \geq \mathcal{J}_n (w) $ and $\mathcal{J}_n$ is sequentially weakly lower semicontinuous. \ \\
It follows that $\mathcal{J}_n$ admits a global minimum at $u_n \in \mathcal{D}_0^{1,p}(\Omega^e)$ for each $n \in \mathbb{N}$. Moreover $\mathcal{J}_n$ is of class $C^1(\mathcal{D}_0^{1,p}(\Omega^e))$ and
\[
\mathcal{J}_n'(u_n)(\varphi) = \int_{\Omega^e }H ^{p-1} \left( \nabla u_n\right) \left\langle \nabla H \left(\nabla u_n \right), \nabla \varphi \right\rangle - \int_{\Omega^e} g_n(x, u_n) \varphi,  \quad \forall \varphi
\in \mathcal{D}_{0}^{1,p}\left( \Omega^e \right) .
\]
Since $\mathcal{J}_n'(u_n)=0$, $u_n$ satisfies \eqref{weak equation appr problem}, so that it is a weak solution of \eqref{appr problem}.
\end{proof}

\begin{lemma}\label{properties of u_n} 
Let $u_{n}\in \mathcal{D}_{0}^{1,p}\left( \Omega^e \right)$ be a weak solution of problem \eqref{appr problem} for each $n\in \mathbb{N}$. Then 
\begin{itemize}
\item[(i)] $u_n\geq 0$ and $u_n(x) \to 0$ as $|x| \to \infty$;
\item[(ii)] $u_n \le u_{n+1}$;
\item[(iii)] $u_n$ is bounded in $W^{1,p}_{loc}(\Omega^e)$ and $\left\{u_n^{\frac{p+\gamma-1}{p}}\right\}_{n \in\N}$ is bounded in $\mathcal{D}_0^{1,p}(\Omega^e)$.
\end{itemize}
\end{lemma}

\begin{proof}
	We first prove (i). Let  $\varphi=u_n^-:=\max\{-u_n,0\} \in \mathcal{D}_0^{1,p}(\Omega^e)$ as test function in \eqref{weak equation appr problem} so we have
$$
 \int_{\Omega^e }H ^{p-1} \left( \nabla u_n\right) \left\langle \nabla H \left(\nabla u_n \right), \nabla u_n^- \right\rangle  
		=\int_{\Omega^e } g_n(x, u_n)u_n^-
$$
We point out that $\nabla u_n^-=-\nabla u_n\,\chi_{\{u_n<0\}}$ and, where ${\{u_n<0\}}$, $g_n(x, u_n)=n^\gamma K(x)$.
Moreover where ${\{u_n<0\}}$
$$
  H^p \left(\nabla u^{-}_n\right) = \left\langle H^{p-1} \left(\nabla u^{-}_n \right) \nabla H \left( \nabla u^{-}_n  \right), \nabla u^{-}_n \right\rangle=-\left\langle H ^{p-1} \left( \nabla u_n\right)  \nabla H \left(\nabla u_n \right), \nabla u_n^- \right\rangle  
$$
hence 
\[
-\int_{\Omega^e}H^p(\nabla u_n^-) = n^\gamma\int_{\Omega^e}K(x)u_n^-\ge0 \implies \int_{\Omega^e}H^p(\nabla u_n^-)=0 .
\]
By \eqref{norm}, $\|u_n^-\|_{\mathcal{D}_0^{1,p}(\Omega^e)}=0$, hence $u_n\geq 0$.
%thus exploiting inequality \eqref{norm} and Euler's theorem \eqref{eq:eulero} we get
%	\begin{align*}
%		a^p|\nabla u^{-}_n|^p \le H^p \left(\nabla u^{-}_n\right) = \left\langle H^{p-1} \left(\nabla u^{-}_n \right) \nabla H \left( \nabla u^{-}_n  \right), \nabla u^{-}_n \right\rangle.
%	\end{align*}
%	Using $\varphi=u^{-}_n$ as test function in \eqref{weak equation appr problem}, by Sobolev inequality we have that
%	\begin{align*}
%		a^{p^*}\int_{\Omega^e } |u^{-}_n|^{p*} \le C_S a^p \int_{\Omega^e } |\nabla u^{-}_n|^p \le C_S \int_{\Omega^e } H^p \left(\nabla u^{-}_n\right) = C_S\int_{\Omega^e} \dfrac{-K(x)}{\left( u^{-}_n+\frac{1}{n}\right)^\gamma}u^{-}_n \le 0.
%	\end{align*}
%%
%Hence, by \eqref{norm}, $a^{p^*}\|u_n^-\|^{p^*}_{L^{p^*}(\Omega^e)}\leq\int_{\Omega^e}H^p(\nabla u_n^-)\leq0$, so that $u_n^-=0$, i.e. $u_n\geq0$. 
For any $n\in\mathbb{N}$ and $\varepsilon\in(0,1)$ fixed, set
\[
u_{\varepsilon,n}\,:=\, (u_n-\varepsilon)^+\,+\,\frac{1}{\varepsilon^\gamma}\,.
\]
Furthermore, for a given $q\geq 1$ and $l >\frac{1}{\varepsilon^\gamma}$, we define 
\begin{equation}
	F(t)\,:=\, \begin{cases}
	t^q\qquad\qquad \qquad\qquad\qquad if\quad \frac{1}{\varepsilon^\gamma}\leq t\leq l\\
	ql^{q-1}t-(q-1)l^q\qquad\quad  if \quad t\geq l
\end{cases}
\end{equation}
and, consequently, 
\begin{equation}
	G(u_n)\,=\, G(u_{\varepsilon,n})\,:=\, F(u_{\varepsilon,n})(F'(u_{\varepsilon,n}))^{p-1}-\frac{q^{p-1}}{\varepsilon^{\beta\gamma}},
\end{equation}
where $\beta:=p(q-1)+1$. We now indicate by $\eta(\cdot)$ a standard smooth cut-off function and consider the test function
\[
\varphi\,:=\,G(u_n)\eta^p\,.
\]
Note that $\varphi$ is a suitable test function in \eqref{weak equation appr problem}. Indeed, by the definitions of $F$ and $\beta$, we have $G(u_n)=0$ in $\{u_n\leq\varepsilon\}$; moreover, $G(u_n)$ is a Lipschitz function of $u_n$, since $F$ is linear for $u_{\varepsilon,n}\geq l$. Hence $\varphi\in\mathcal{D}_0^{1,p}(\Omega^e)$ and it has bounded support. Consequently we deduce that
\begin{equation}\nonumber
\begin{split}
\int_{\Omega^e}\, H^p(\nabla u_n) G'(u_n) \eta^p&=-p\int_{\Omega^e}G(u_n)\eta^{p-1} H^{p-1}(\nabla u_n)\langle\nabla H(\nabla u_n)\,,\nabla\eta\rangle\\
&\quad+\int_{\Omega^e}g_n(x,u_n)G(u_n)\eta^p%\frac{K(x)G(u_n)\eta^p}{u_n^\gamma}
\end{split}
\end{equation}
so that 
\begin{equation}\nonumber
	\begin{split}
		\int_{\Omega^e}\, H^p(\nabla u_n) G'(u_n)\eta^p&\leq p \int_{\Omega^e} H^{p-1}(\nabla u_n)G(u_n)\eta^{p-1} |\nabla H(\nabla u_n)||\nabla\eta|\\
		&\quad+\int_{\Omega^e}\frac{K(x)G(u_n)\eta^p}{u_n^\gamma}\\
		&\leq p \int_{\Omega^e} H^{p-1}(\nabla u_n)G(u_n)\eta^{p-1}|\nabla H(\nabla u_n)||\nabla\eta|\\
		&\quad+ \int_{\Omega^e}\frac{K(x)G(u_n)\eta^p}{\varepsilon^{\gamma}} \\
		&\leq p \int_{\Omega^e} H^{p-1}(\nabla u_n)G(u_n)\eta^{p-1}|\nabla H(\nabla u_n)||\nabla\eta|\\
		&\quad+ c_\varepsilon\int_{\Omega^e}K(x)G(u_n)\eta^pu_{\varepsilon,n}^{p-1}, \\
	\end{split}
\end{equation}
where $c_\varepsilon:=\max\{1,\varepsilon^{-\gamma(2-p)}\}$ and we have used that $g_n(x,u_n)\leq K(x)u_n^{-\gamma}$, that $G(u_n)=0$ in $\{u_n\leq\varepsilon\}$ and that, in $\{u_n>\varepsilon\}$, we have $u_n^{-\gamma}<\varepsilon^{-\gamma}$ and
$$\varepsilon^{-\gamma}\leq c_\varepsilon\,\varepsilon^{-\gamma(p-1)}\leq c_\varepsilon\left((u_n-\varepsilon)^+\,+\,\frac{1}{\varepsilon^\gamma}\right)^{p-1}=c_\varepsilon\,u_{\varepsilon,n}^{p-1}.$$
Taking into account the definition of $F(\cdot)$ and $G(\cdot)$, it is easily seen that
	\begin{equation}\nonumber
	G'(u_n)\,=\, \begin{cases}
		\dfrac{\beta}{q}\,(F'(u_{\varepsilon, n}))^p & \text{if }\ 0<u_n-\varepsilon<l-\frac{1}{\varepsilon^\gamma},\\[3mm]
		(F'(u_{\varepsilon, n}))^p & \text{if }\ u_n-\varepsilon>l-\frac{1}{\varepsilon^\gamma},\\[2mm]
		0 & \text{if }\ u_n<\varepsilon,
	\end{cases}
\end{equation}
and $G\leq F(F')^{p-1}$. Using this last we can see that
\begin{equation}\nonumber
	\begin{split}
		\min \left\{1,\frac{\beta}{q}\right\}\int_{\Omega^e}\, (H(\nabla u_n) F'(u_{\varepsilon, n})\eta)^p&\leq p \int_{\Omega^e} (H(\nabla u_n)F'(u_{\varepsilon, n}))^{p-1}F(u_{\varepsilon, n})\eta^{p-1} |\nabla H(\nabla u_n)||\nabla\eta|\\
		&\quad + c_\varepsilon \int_{\Omega^e} K(x) F(u_{\varepsilon, n})(F'(u_{\varepsilon, n}) u_{\varepsilon, n})^{p-1}\eta^p\,. \\
	\end{split}
\end{equation}
Thus, setting
\[
v(x)\,:=\,F(u_{\varepsilon,n}) \qquad (\Rightarrow \nabla v=F'(u_{\varepsilon,n}) \nabla u_{\varepsilon,n})
\]
and observing that $u_{\varepsilon,n}\,F'\leq q F$,  by \eqref{norm} and \eqref{grad H limitato}, we get that
\begin{equation}\nonumber
	\begin{split}
	\min \left\{a^p,\frac{a^p\beta}{q}\right\}\int_{\Omega^e}|\nabla v|^p\eta^p
		&= a^p\min \left\{1,\frac{\beta}{q}\right\}\int_{\Omega^e}(|\nabla u_{\varepsilon,n}| F'(u_{\varepsilon,n})\eta)^p\\
		&\leq \min \left\{1,\frac{\beta}{q}\right\}\int_{\Omega^e}(H(\nabla u_{\varepsilon,n})F'(u_{\varepsilon,n})\eta)^p\\
		& \leq pb^{p-1}C_H \int_{\Omega^e} (|\nabla v|\eta)^{p-1}|\nabla\eta|v\\
		&\quad +c_\varepsilon q^{p-1}\int_{\Omega^e}K(x)v^p\eta^p\,\\
		%& \leq pb^{p-1}C_H \|v\nabla \eta\|_p\|\eta \nabla v\|^{p-1}_p+q^{p-1}\int_{\Omega^e}K(x)v^p\eta^p\,.
	\end{split}
\end{equation}
%By a simple application of the Young inequality and then by Sobolev inequality, we observe that
%$$ 
%\int_{\Omega^e}K(x)v^p\eta^p\leq \|K\|_{\frac{N}{p-\epsilon}}\|\eta v\|_ p^\epsilon\|\eta v\|_{p^*}^{p-\epsilon}\leq C_S\|K\|_{\frac{N}{p-\epsilon}}\|\eta v\|_ p^\epsilon\|\nabla(\eta v)\|_{p}^{p-\epsilon}.
%$$
%Moreover recall that there exists $C_p>0$ such that
%$$
% |\nabla (v\eta)|^p\leq C_p|\nabla \eta|^pv^p+C_p|\nabla v|^p\eta^p, 
%$$
%therefore, redefining the constant
%\begin{equation}\nonumber
%	\begin{split}
%	\min \left\{a^p,\frac{a^p\beta}{q}\right\}\int_{\Omega^e}|\nabla v|^p\eta^p&\leq C_pb^{p-1}C_H \|v\nabla \eta\|_p\|\eta \nabla v\|^{p-1}_p\\
%	&\quad +C_{p,\epsilon}C_Sq^{p-1}\|K\|_{\frac{N}{p-\epsilon}}\|\eta v\|_ p^\epsilon (\|v\nabla \eta\|_p^{p-\epsilon}+\|\eta\nabla v\|_p^{p-\epsilon}).
%	\end{split}
%\end{equation}

This estimate recovers equation (16) in the proof of Theorem 1 of \cite{Se} if $\bar{b}=c=0$. Consequently the remaining part of the proof can be carried over repeating verbatim the argument of Serrin, thus proving the local boundedness of $u_{\varepsilon,n}$ that implies the local boundedness of $u_n$.

Using standard regularity results (recall that, by \eqref{hyp:K}, $0<g_n(\cdot,u_n)\leq n^\gamma K\in L^q(\Omega^e)$ for every $q\in[1,\infty)$), we have $u_n \in C^{1,\alpha}(\Omega^e)$ locally, i.e. for each $r > 0$ there exists $\alpha \in (0,1)$ such that $u_n \in C^{1,\alpha}(\Omega^e \cap B_r)$. By the strong maximum principle we get $u_n > 0$ in $\Omega^e$ for any $n\in\mathbb{N}$. Moreover, for $|x|$ large enough, so that $B_2(x)\subset\Omega^e$, the local $L^{\infty}$-estimates (see \cite[Theorem 1]{Se}) give
\begin{equation*}
	\|u_n\|_{L^\infty(B_1(x))} \le C\Big(\|u_n\|_{L^{p}(B_2(x))}+\big(n^\gamma\|K\|_{L^{\frac{N}{p-\delta}}(B_2(x))}\big)^{\frac{1}{p-1}}\Big),
\end{equation*}

%where $\delta\in(0,1)$ is fixed and $C=C(p,N,a,b,C_H,\delta)$ does not depend on $x$. Since $u_n\in L^{p^*}(\Omega^e)$ and $K\in L^{\frac{N}{p-\delta}}(\Omega^e)$, the right-hand side tends to $0$ as $|x|\to\infty$; hence
%
\begin{equation}\label{eq:decay-un}
	\lim_{|x|\to\infty} u_n(x) = 0 \quad \text{for each } n \in \mathbb{N}.
\end{equation}

Now we show \emph{(ii)} that is $u_n \le u_{n+1}$. Using $\varphi = (u_{n+1} - u_n)^-$ as test function in \eqref{weak equation appr problem}, we have that
\begin{equation*}
	\begin{split}
	\int_{\Omega^e }& \left\langle H ^{p-1} \left( \nabla u_{n}\right)\nabla H \left(\nabla u_{n} \right)-H ^{p-1} \left( \nabla u_{n+1}\right)\nabla H \left(\nabla u_{n+1} \right), \nabla \left( u_{n+1}-u_n\right)^-  \right\rangle \\
	&=\int_{\Omega^e} K(x) \left[  \left(u_n + \frac{1}{n}\right)^{-\gamma}-\left(u_{n+1} + \frac{1}{n+1}\right)^{-\gamma}  \right] (u_{n+1} - u_n)^-.
	\end{split}
\end{equation*}
By \eqref{1.4}-Lemma \ref{estimates} we get
\begin{equation*}
\begin{split}
	C_1\int_{\Omega^e } & \left( |\nabla u_n|+|\nabla u_{n+1}|\right)^{p-2} |\nabla(u_{n+1}-u_n)^-|^2 \\
	&\leq \int_{\Omega^e} K(x) \left[ \left(u_n + \frac{1}{n}\right)^{-\gamma} - \left(u_{n+1} + \frac{1}{n+1}\right)^{-\gamma} \right] (u_{n+1} - u_n)^-.
\end{split}
\end{equation*}
It follows
\begin{equation*}
\begin{split}
	0 &\le C_1\int_{\text{supp}( (u_{n+1} - u_n)^-)}\left( |\nabla u_n|+|\nabla u_{n+1}|\right)^{p-2} |\nabla(u_n-u_{n+1})|^2 \\
	&\leq\int_{\text{supp}( (u_{n+1} - u_n)^-)} K(x) \left[ \left(u_n + \frac{1}{n}\right)^{-\gamma}- \left(u_{n+1} + \frac{1}{n+1}\right)^{-\gamma} \right] (u_n-u_{n+1}) \le 0,
\end{split}
\end{equation*}
hence $u_n \le u_{n+1}$ for each $n\in \mathbb{N}$. \\
Then, since $u_1 \in C(\Omega^e)$ and $u_1 > 0$ in $\Omega^e$, for each compact set $\mathcal{K} \subset \Omega^e$ there is a constant $c_\mathcal{K}> 0$ such that $u_n \ge u_1 > c_\mathcal{K}$ in $\mathcal{K}$ for each $n \in \mathbb{N}$.\\ \\
Finally, we prove \emph{(iii)}, i.e. $\{u_n\}$ is bounded in $W_{loc}^{1,p}(\Omega^e)$. We start showing that $\left\{u_n^{\frac{p+\gamma-1}{p}}\right\}_{n \in\mathbb{N}}$ is bounded in $\mathcal D_0^{1,p}(\Omega^e)$.  Using \(T_M(u_n)^\gamma\) as test function\footnote{$
T_M(s)=\min\{s,M\}$}  in \eqref{weak equation appr problem} and then letting \(M\to\infty\), we obtain,  by \eqref{norm}, 
\begin{align*}
	\int_{\Omega^e} \left|\nabla \left(u_n^{\frac{p+\gamma-1}{p}}\right)\right|^p dx &\leq \left( \frac{p+\gamma-1}{p} \right)^p \int_{\Omega^e} |\nabla u_n|^p u_n^{\gamma-1} dx \\ \mbox{(by \eqref{norm})}&\le \left( \frac{p+\gamma-1}{ap} \right)^p \int_{\Omega^e} H^p \left(\nabla u_n\right)u_n^{\gamma-1} \\
	&= \frac{1}{\gamma} \left( \frac{p+\gamma-1}{ap} \right)^p \int_{\Omega^e}H ^{p-1} \left( \nabla u_n\right) \left\langle \nabla H \left(\nabla u_n \right), \nabla u_n^\gamma \right\rangle \\ &= \frac{1}{\gamma} \left( \frac{p+\gamma-1}{ap} \right)^p \int_{\Omega^e} g_n(x, u_n) u_n^\gamma dx \leq  \frac{ \|K\|_1}{\gamma} \left( \frac{p+\gamma-1}{ap} \right)^p,
\end{align*}
where we took into account the definition of $g_n$. This implies that $\left\{u_n^{\frac{p+\gamma-1}{p}}\right\}$ is bounded in $\mathcal{D}^{1,p}_0(\Omega^e)$.\\
Moreover, if $\mathcal{K} \subset \Omega^e$ is a compact set and $c_\mathcal{K}  > 0$ is such that $u_n > c_\mathcal{K} $ in $\mathcal{K} $ for each $n \in \mathbb{N}$, we deduce that
\begin{equation*}
	c_\mathcal{K} ^{\gamma-1} \int_\mathcal{K} |\nabla u_n|^p dx \le \int_{\Omega^e} u_n^{\gamma-1} |\nabla u_n|^p dx \le \frac{\|K\|_1}{a\gamma}.% \left( \frac{p}{p+\gamma-1} \right)^p.
\end{equation*}
This implies that $\{u_n\}$ is bounded in $W_{loc}^{1,p}(\Omega^e)$.
 \end{proof}

\begin{lemma}\label{lem:sandwich}
Let $\gamma>0$ and, for every $n\in\mathbb N^+$, let $u_n\in\mathcal D^{1,p}_0(\Omega^e)$ be a weak solution of \eqref{appr problem}, i.e.\ let $u_n$ satisfy \eqref{weak equation appr problem}. Then
\begin{equation}\label{eq:sandwich-mn}
	u_m\le u_n+\frac1n\quad\text{a.e. in }\Omega^e,\qquad\text{for every }m,n\in\mathbb N^+ .
\end{equation}
Consequently, the function $u:=\lim_{n\to\infty}u_n=\sup_{n}u_n$ (see Lemma~\ref{properties of u_n}-(ii)) satisfies
\begin{equation}\label{eq:sandwich}
	u_n\le u\le u_n+\frac1n\quad\text{a.e. in }\Omega^e,\qquad\text{for every }n\in\mathbb N^+ .
\end{equation}
In particular, $u_n\to u$ uniformly in $\Omega^e$, $u\in C(\Omega^e)$   and $\lim_{|x|\to\infty}u(x)=0$.
\end{lemma}
\begin{proof}
		Fix $m,n\in\mathbb N^+$ and set
		\[
		\varphi:=\Big(u_m-u_n-\frac1n\Big)^+ .
		\]
		Since $u_m-u_n\in\mathcal D^{1,p}_0(\Omega^e)$ and $\frac1n>0$, we have
		$\varphi\in\mathcal D^{1,p}_0(\Omega^e)$, with
		$\nabla\varphi=\chi_{\{u_m-u_n>1/n\}}\nabla(u_m-u_n)$. Since $u_m,u_n\ge0$
		(Lemma~\ref{properties of u_n}-(i)), we have
		$g_m(x,u_m)=K(x)\big(u_m+\frac1m\big)^{-\gamma}$ and
		$g_n(x,u_n)=K(x)\big(u_n+\frac1n\big)^{-\gamma}$. Using $\varphi$ as test
		function in \eqref{weak equation appr problem}, written both for $u_m$ and for $u_n$, and
		subtracting, we get
		\begin{multline}\label{eq:sandwich-proof}
			\int_{\Omega^e \cap \text{supp}(\varphi)}\Big\langle H^{p-1}(\nabla u_m)\nabla H(\nabla u_m)
			-H^{p-1}(\nabla u_n)\nabla H(\nabla u_n),\nabla(u_m-u_n)\Big\rangle\\
			=\int_{\Omega^e}K(x)\left[\frac{1}{\big(u_m+\frac1m\big)^{\gamma}}
			-\frac{1}{\big(u_n+\frac1n\big)^{\gamma}}\right]\Big(u_m-u_n-\frac1n\Big)^+,
		\end{multline}
		where both sides are finite, since $\varphi\in L^{p^*}(\Omega^e)$ and
		$K\in L^{(p^*)'}(\Omega^e)$. On the set $\{\varphi>0\}$ we have
		$u_m+\frac1m>u_m>u_n+\frac1n$, so that the right-hand side of
		\eqref{eq:sandwich-proof} is non-positive. On the other hand, by
		\eqref{1.4}, we get
		\begin{multline}\nonumber
			0 \leq C_1\int_{\Omega^e\cap \text{supp}(\varphi)}\big(|\nabla u_m|+|\nabla u_n|\big)^{p-2}|\nabla(u_m-u_n)|^2\  \\
			\leq \int_{\Omega^e \cap \text{supp}(\varphi)}\Big\langle H^{p-1}(\nabla u_m)\nabla H(\nabla u_m)
			-H^{p-1}(\nabla u_n)\nabla H(\nabla u_n),\nabla(u_m-u_n)\Big\rangle \leq 0.
		\end{multline}
		Hence 
		\[
		C_1\int_{\Omega^e}\big(|\nabla u_m|+|\nabla u_n|\big)^{p-2}|\nabla\varphi|^2\ =0,
		\]
		and, since $\nabla\varphi=0$ a.e.\ on the
		set where $|\nabla u_m|+|\nabla u_n|=0$, we infer that $\nabla\varphi=0$
		a.e.\ in $\Omega^e$, i.e.\ $\|\varphi\|=0$. Thus $\varphi=0$, and
		\eqref{eq:sandwich-mn} is proved.
		Letting $m\to\infty$ in \eqref{eq:sandwich-mn} we obtain the second
		inequality in \eqref{eq:sandwich}, while the first one follows from the
		monotonicity of $(u_n)$. Therefore $\|u-u_n\|_{L^\infty(\Omega^e)}\le\frac1n$;
		since each $u_n$ is continuous in $\Omega^e$, we deduce that $u_n\to u$ uniformly in $\Omega^e$
		and that $u$ is continuous in $\Omega^e$. Finally, fix $\varepsilon>0$ and
		$\bar n\in\mathbb N^+$ such that $\frac1{\bar n}<\frac\varepsilon2$. By
		\eqref{eq:decay-un} there exists $\rho>0$ such that
		$u_{\bar n}(x)<\frac\varepsilon2$ for $|x|\ge\rho$; then, by
		\eqref{eq:sandwich},
		\[
		0\le u(x)\le u_{\bar n}(x)+\frac1{\bar n}<\varepsilon\qquad\text{for }|x|\ge\rho ,
		\]
		that is, $\lim_{|x|\to\infty}u(x)=0$.
\end{proof}

\begin{theorem}\label{existence}
	Let $\gamma \geq 1$ and let $K$ be a positive function satisfying \eqref{hyp:K}. Then problem \eqref{problem} admits a solution $u$, in the sense of Definition~\ref{def:sol}, such that $u^{\frac{p+\gamma-1}{p}}\in\mathcal{D}_0^{1,p}(\Omega^e)$. 
\end{theorem}
\begin{proof}
	By Lemma \ref{properties of u_n} there exists $u \in W_{loc}^{1,p}(\Omega^e)$ such that, up to a subsequence, $u_n \rightharpoonup u$ in $W^{1,p}_{loc}(\Omega^e)$; since $(u_n)$ is nondecreasing, $u=\sup_n u_n$ a.e. in $\Omega^e$. Let us prove that $u$ is a solution of problem \eqref{problem}. To this aim let us fix $\varphi \in C_c^1(\Omega^e)$ and let $c_\varphi > 0$ be such that $u_n\geq u_1 \ge c_\varphi$ on $\text{supp}(\varphi)$. Since
	\begin{equation*}
		g_n(x, u_n(x))\varphi(x) \to \frac{K(x)}{u(x)^\gamma}\varphi(x) \text{ for a.e. } x \in \Omega^e
	\end{equation*}
	and
	\begin{equation*}
			|g_n(x, u_n(x))\varphi(x)| \le \frac{K(x)}{c_\varphi^\gamma}|\varphi(x)| \in L^1(\Omega^e),
	\end{equation*}
	by the Lebesgue's dominated convergence theorem we have $\int_{\Omega^e}g_n(x,u_n)\varphi\to\int_{\Omega^e}\frac{K(x)}{u^\gamma}\varphi$. On the other hand, $(u_n)$ is bounded in $W^{1,p}_{loc}(\Omega^e)$ and, for every compact set $\omega\subset\Omega^e$, we have $0<g_n(x,u_n)\leq c_\omega^{-\gamma}K(x)$ in $\omega$, so that $(g_n(\cdot,u_n))$ is bounded in $L^1_{loc}(\Omega^e)$. Hence, by \cite[Theorem 2.1 and Remark 2.2]{BM} (see also \cite{DMM,EST}), we have $\nabla u_n\to\nabla u$ a.e. in $\Omega^e$. Since, by \eqref{norm} and \eqref{grad H limitato}, $\big(H^{p-1}(\nabla u_n)\nabla H(\nabla u_n)\big)$ is bounded in $L^{p'}_{loc}(\Omega^e;\mathbb{R}^N)$, the Vitali's theorem yields
	\[
	\int_{\Omega^e }H ^{p-1} \left( \nabla u_n\right) \left\langle \nabla H \left(\nabla u_n \right), \nabla \varphi \right\rangle\to\int_{\Omega^e }H ^{p-1} \left( \nabla u\right) \left\langle \nabla H \left(\nabla u \right), \nabla \varphi \right\rangle .
	\]
	Therefore, passing to the limit in
	\begin{equation*}
	\int_{\Omega^e }H ^{p-1} \left( \nabla u_n\right) \left\langle \nabla H \left(\nabla u_n \right), \nabla \varphi \right\rangle   = \int_{\Omega^e} g_n(x, u_n)\varphi  \quad \text{for all } n \in \mathbb{N}^+,
	\end{equation*}
	we obtain
	\begin{equation*}
			\int_{\Omega^e }H ^{p-1} \left( \nabla u\right) \left\langle \nabla H \left(\nabla u \right), \nabla \varphi \right\rangle   = \int_{\Omega^e} \frac{K(x)}{u(x)^\gamma}\varphi.
	\end{equation*}
Moreover, since $u$ is the pointwise limit of $\{u_n\}$, for each compact set $\omega \subset \Omega^e$, there exists $c_\omega > 0$ such that $u > c_\omega$ in $\omega$. Furthermore, by Lemma~\ref{properties of u_n}-(iii), $\left\{u_n^{\frac{p+\gamma-1}{p}}\right\}_{n \in \N}$ is bounded in $\mathcal{D}_0^{1,p}(\Omega^e)$; since it converges a.e. to $u^{\frac{p+\gamma-1}{p}}$, we infer that $u^{\frac{p+\gamma-1}{p}}\in\mathcal{D}_0^{1,p}(\Omega^e)$. Now, let $\{\varphi_n\} \subset C_c^1(\Omega^e)$ be a sequence of nonnegative functions converging to $u^{\frac{p+\gamma-1}{p}}$ in $\mathcal{D}_0^{1,p}(\Omega^e)$. Then, the sequence $\psi_n := (\varphi_n^{\frac{p}{p+\gamma-1}} - \varepsilon)^+$ converges a.e. to $(u - \varepsilon)^+$. Moreover,
	\begin{align*}
		\int_{\Omega^e} |\nabla \psi_n|^p dx &= \left( \frac{p}{p+\gamma-1} \right)^p \int_{\varphi_n^{\frac{p}{p+\gamma-1}} > \varepsilon} (\varphi_n^{\frac{p}{p+\gamma-1}})^{1-\gamma} |\nabla \varphi_n|^p dx \\
		&\le \left( \frac{p}{p+\gamma-1} \right)^p \varepsilon^{1-\gamma} \int_{\Omega^e} |\nabla \varphi_n|^p dx < C,
	\end{align*}
	for some positive constant $C$. Then, since $\{\psi_n\}$ is uniformly bounded in $\mathcal{D}_0^{1,p}(\Omega^e)$ and converges a.e. to $(u - \varepsilon)^+$, we have that $(u - \varepsilon)^+ \in \mathcal{D}_0^{1,p}(\Omega^e)$. Fix $\varepsilon > 0$ and $\bar{n}$ such that $\frac{1}{\bar{n}} < \frac{\varepsilon}{2}$. Using again Lemma \ref{properties of u_n} we can find $\rho > 0$ such that $u_{\bar{n}}(x) < \frac{\varepsilon}{2}$ for $|x| \ge \rho$. By Lemma~\ref{lem:sandwich} we have
	\begin{equation*}
		u_n \le u \le u_n + \frac{1}{n} \quad \text{in } \Omega^e \quad \text{for all } n \in \mathbb{N}^+,
	\end{equation*}	
	hence $0 \le u(x) \le u_{\bar{n}}(x) + \frac{1}{\bar{n}} < \varepsilon$ for $|x| \ge \rho$, so that $\lim_{|x|\to\infty} u(x) = 0$. Thus, the proof is complete.
\end{proof}

\subsection{The case $0<\gamma<1$.} In what follows let us introduce 
$$
 m:=\frac{Np}{N\left( p-1\right) +p+\gamma \left( N-p\right) }=\left( \frac{%
			p^{\ast }}{1-\gamma }\right) ^{\prime },
$$
Let us immediately remark that 
\begin{itemize}
 \item $\displaystyle 1<m<\frac {N}{p}<\frac {N}{p^-}$, hence, by interpolation, $K\in L^m(\Omega^e)$.\ \\
 \item $\displaystyle m':=\frac{m}{m-1}$ is such that $(1-\gamma)m'=p^*$.
\end{itemize}

\noindent The proof of the following lemma is the same as that of \cite[Lemma 2.7]{EST}, since it only relies on the H\"older and Sobolev inequalities.

\begin{lemma}\label{gammasotto1unif}
	Let $0<\gamma<1$. If $u_{n}$ solves \eqref{appr problem} and $K$ is a positive function satisfying \eqref{hyp:K}, then $\left( u_{n}\right) $ is uniformly bounded in $\mathcal{D}_{0}^{1,p}\left(
	\Omega^e \right) $.
\end{lemma}

\begin{theorem} \label{existence 3}
		Let $0<\gamma<1$. If $K$ is a positive function satisfying \eqref{hyp:K}, then \eqref{problem} has a solution $u\in \mathcal{D}_{0}^{1,p}\left( \Omega^e\right)$, in the sense of Definition~\ref{def:sol}.
\end{theorem}

\begin{proof}
		Let $(u_n)$ be the sequence of solutions of \eqref{appr problem} given by
		Lemma~\ref{existence appr problem}.	By \eqref{norm} and \eqref{grad H limitato},
		$H^{p-1}(\xi)|\nabla H(\xi)| \le b^{p-1}C_H|\xi|^{p-1}$ for every $\xi\in\mathbb R^N$, so
		that $H^{p-1}(v)|\nabla H(v)|\in L^{p'}(\Omega^e)$ for every
		$v\in\mathcal D^{1,p}_0(\Omega^e)$. We also point out that items (i)--(ii) of
		Lemma~\ref{properties of u_n}, as well as the lower bound on compact sets obtained
		in its proof, do not depend on $\gamma$. Hence $0\le u_n\le u_{n+1}$ in
		$\Omega^e$ for every $n$, and for every compact set $\omega\subset\Omega^e$
		there exists $c_\omega>0$ such that
		\begin{equation}\label{eq:lower-bound}
			u_n\ge u_1\ge c_\omega\quad\text{in }\omega,\qquad\text{for every }n\in\mathbb N^+ .
		\end{equation}
		By Lemma~\ref{gammasotto1unif} there
		exists $C>0$ such that $\|u_n\|\le C$ for every $n$. Let
		$u:=\lim_{n}u_n=\sup_n u_n$, which exists by monotonicity. By the Sobolev inequality and the monotone convergence theorem,
		\[
		\|u\|_{p^*}=\lim_{n\to\infty}\|u_n\|_{p^*}\le C_S\,C ,
		\]
		so that $u$ is finite a.e.\ in $\Omega^e$. Moreover, since
		$0\le u-u_n\le u\in L^{p^*}(\Omega^e)$, the dominated convergence theorem
		yields
		\begin{equation}\label{eq:strong-Lpstar}
			u_n\to u\quad\text{in }L^{p^*}(\Omega^e).
		\end{equation}
		Since $\mathcal D^{1,p}_0(\Omega^e)$ is reflexive and continuously embedded in
		$L^{p^*}(\Omega^e)$, every subsequence of $(u_n)$ has a further subsequence
		which converges weakly in $\mathcal D^{1,p}_0(\Omega^e)$, and by
		\eqref{eq:strong-Lpstar} its weak limit is $u$. Therefore
		$u\in\mathcal D^{1,p}_0(\Omega^e)$, the whole sequence satisfies
		$u_n\rightharpoonup u$ in $\mathcal D^{1,p}_0(\Omega^e)$, and
		$\|u\|\le\liminf_n\|u_n\|\le C$.
		
		\
		
		Since $u\in\mathcal D^{1,p}_0(\Omega^e)$, the function $u_n-u$ is an admissible test
		function in \eqref{weak equation appr problem}. As $u_n\le u$ and $g_n>0$, we get
		\[
		\int_{\Omega^e} H^{p-1}(\nabla u_n) \langle  \nabla H(\nabla u_n), \nabla(u_n-u)\rangle =\int_{\Omega^e}g_n(x,u_n)\,(u_n-u)\le0 ,
		\]
		where the right-hand side is finite, since
		$0<g_n(x,u_n)\le n^\gamma K\in L^{(p^*)'}(\Omega^e)$. Then, by
		\eqref{1.4},
		\begin{multline*}
			0\le C_1\int_{\Omega^e}\big(|\nabla u_n|+|\nabla u|\big)^{p-2}|\nabla(u_n-u)|^2\\
			\le\int_{\Omega^e}\langle H^{p-1}(\nabla u_n) \nabla H(\nabla u_n)-H^{p-1}(\nabla u) \nabla H(\nabla u),\nabla(u_n-u)\rangle\\
			\le-\int_{\Omega^e}\langle  H^{p-1}(\nabla u) \nabla H(\nabla u),\nabla(u_n-u)\rangle ,
		\end{multline*}
		and the last term tends to $0$, because $\nabla u_n\rightharpoonup\nabla u$
		in $L^p(\Omega^e;\mathbb R^N)$ and
		$H^{p-1}(\nabla u) |\nabla H(\nabla u)| \in L^{p'}(\Omega^e)$. Hence
		\begin{equation}\label{eq:In}
			I_n:=\int_{\Omega^e}\big(|\nabla u_n|+|\nabla u|\big)^{p-2}|\nabla(u_n-u)|^2
			\longrightarrow0 .
		\end{equation}
		If $p\ge2$, since $|\nabla(u_n-u)|\le|\nabla u_n|+|\nabla u|$, we have
		$\|u_n-u\|^p\le I_n$. If $1<p<2$, the H\"older inequality with exponents
		$\frac2p$ and $\frac2{2-p}$ gives
		\[
		\|u_n-u\|^p\le I_n^{\frac p2}\left(\int_{\Omega^e}\big(|\nabla u_n|+|\nabla u|\big)^p\right)^{\frac{2-p}2}
		\le(2C)^{\frac{p(2-p)}2}\,I_n^{\frac p2} .
		\]
		In both cases, \eqref{eq:In} implies that $u_n\to u$ strongly in
		$\mathcal D^{1,p}_0(\Omega^e)$.
		
		Finally, we prove that $u$ satisfies \eqref{weak equation}. Fix
		$\varphi\in C^1_c(\Omega^e)$ and let $\omega:=\operatorname{supp}\varphi$.
		Since $\varphi\in\mathcal D^{1,p}_0(\Omega^e)$, by \eqref{weak equation appr problem} we
		have
		\begin{equation}\label{eq:weak-n-phi}
			\int_{\Omega^e}\langle H^{p-1}(\nabla u_n) \nabla H(\nabla u_n),\nabla\varphi\rangle
			=\int_{\Omega^e}\frac{K(x)}{\big(u_n+\frac1n\big)^{\gamma}}\,\varphi
			\qquad\text{for every }n\in\mathbb N^+ .
		\end{equation}
		Let us pass to the limit in the left-hand side. By \eqref{1.5},
		\begin{multline*}
		\left|\int_{\Omega^e}\langle H^{p-1}(\nabla u_n) \nabla H(\nabla u_n)-H^{p-1}(\nabla u) \nabla H(\nabla u),\nabla\varphi\rangle\right|\\
		\le C_2\int_{\Omega^e}\big(|\nabla u_n|+|\nabla u|\big)^{p-2}|\nabla(u_n-u)|\,|\nabla\varphi| .
		\end{multline*}
		If $p>2$, we apply the H\"older inequality with exponents $\frac p{p-2}$,
		$p$, $p$ (if $p=2$ it is straighfoward). If $1<p<2$, since $|\nabla(u_n-u)|\le|\nabla u_n|+|\nabla u|$ we
		have
		\[
		\big(|\nabla u_n|+|\nabla u|\big)^{p-2}|\nabla(u_n-u)|\le|\nabla(u_n-u)|^{p-1},
		\]
		and we apply the H\"older inequality with exponents $p'$ and $p$. Thus the
		right-hand side above is bounded by
		\[
		\begin{cases}
			C_2(2C)^{p-2}\,\|u_n-u\|\,\|\nabla\varphi\|_p & \text{if }p\ge2,\\[2pt]
			C_2\,\|u_n-u\|^{p-1}\,\|\nabla\varphi\|_p & \text{if }1<p<2,
		\end{cases}
		\]
		and therefore, since $u_n\to u$ strongly in
		$\mathcal D^{1,p}_0(\Omega^e)$, we deduce that
		\[
		\int_{\Omega^e} H^{p-1}(\nabla u_n) \langle \nabla H(\nabla u_n),\nabla\varphi\rangle\longrightarrow
		\int_{\Omega^e}H^{p-1}(\nabla u) \langle \nabla H(\nabla u),\nabla\varphi\rangle .
		\]
		As for the right-hand side of \eqref{eq:weak-n-phi}, by
		\eqref{eq:lower-bound} we have
		\[
		\left|\frac{K(x)}{\big(u_n+\frac1n\big)^{\gamma}}\,\varphi\right|
		\le\frac{\|\varphi\|_{\infty}}{c_\omega^{\gamma}}\,K\,\chi_\omega\in L^1(\Omega^e),
		\]
		and $K\big(u_n+\frac1n\big)^{-\gamma}\varphi\to Ku^{-\gamma}\varphi$ a.e.\ in
		$\Omega^e$, because $u_n\to u$ a.e.\ and $u\ge c_\omega>0$ in $\omega$. By the
		Lebesgue dominated convergence theorem, passing to the limit in \eqref{eq:weak-n-phi} we obtain
		\[
		\int_{\Omega^e}H^{p-1}(\nabla u)\langle\nabla H(\nabla u),\nabla\varphi\rangle
		=\int_{\Omega^e}\frac{K(x)}{u^{\gamma}}\,\varphi
		\qquad\text{for every }\varphi\in C^1_c(\Omega^e),
		\]
		that is, (i) of Definition~\ref{def:sol} holds.
		
		\
		
		Condition (ii) of Definition~\ref{def:sol}  follows from \eqref{eq:lower-bound}, since $u\ge u_n$ for
		every $n$. Condition (iii) holds because $u\ge0$ and
		$u\in\mathcal D^{1,p}_0(\Omega^e)$, so that
		$(u-\varepsilon)^+\in\mathcal D^{1,p}_0(\Omega^e)$ for every $\varepsilon>0$.
		Finally, (iv) follows from Lemma~\ref{lem:sandwich}.
		Hence $u\in\mathcal D^{1,p}_0(\Omega^e)$ is a solution of \eqref{problem} in
		the sense of Definition~\ref{def:sol}.
	\end{proof}

\begin{remark}
		The proof actually shows that $u_n\to u$ both strongly in
		$\mathcal D^{1,p}_0(\Omega^e)$ and uniformly in $\Omega^e$; in particular
		$u\in C(\Omega^e)$.
\end{remark}

\begin{proof}[Proof of Theorem~\ref{thm:existence}]
It follows from Theorems~\ref{existence} and \ref{existence 3}.
\end{proof}

\section{Uniqueness}\label{sec:uniqueness}

The last section of the present paper deals with uniqueness of solutions to problem \eqref{problem}. We recall here the notion of sub and supersolution.

\begin{definition}\label{def:subsuper}
We say that $u\in W^{1,p}_{loc}(\Omega^e)$ is a \emph{subsolution} (respectively, a \emph{supersolution}) of \eqref{problem} if $u$ satisfies (ii) of Definition~\ref{def:sol} and
\[
\int_{\Omega^e}H^{p-1}(\nabla u)\langle\nabla H(\nabla u),\nabla\varphi\rangle\leq\int_{\Omega^e}\frac{K(x)}{u^{\gamma}}\,\varphi\qquad\Big(\text{respectively, }\geq\Big)
\]
for every $\varphi\in C^1_c(\Omega^e)$ with $\varphi\geq0$; moreover, a subsolution is also required to satisfy (iii) of Definition~\ref{def:sol}. In particular, every solution of \eqref{problem} in the sense of Definition~\ref{def:sol} is both a subsolution and a supersolution of \eqref{problem}.
\end{definition}

\noindent Let us define the following real valued function 
\begin{equation}\nonumber
	g_{k}(s)\,:=\,
	\begin{cases}
		\max\{-s^{-\gamma}\,,\,-k\}\quad\text{if}\quad s>0,\\
		-k\qquad\qquad\quad\qquad\text{if}\quad s\leqslant 0\,.\\
	\end{cases}
\end{equation}
We also consider the primitive of $g_{k}(s)$ that is equal to zero for $s=1$, namely
\begin{equation}\nonumber
	\begin{cases}
		\Phi_{k}'(s)\,=\,g_{k}(s),\\
		\Phi_{k}(1)=0\,,\\
	\end{cases}
\end{equation}
and the functional $\mathcal{J}_k\,:\mathcal{D}_0^{1,p}(\Omega^e)\rightarrow [-\infty,+\infty]$ defined by
\begin{equation}\label{Jk}
	\mathcal{J}_k(\varphi)\,:=\,\frac{1}{p}\int_{\Omega^e}H ^{p} \left( \nabla \varphi\right)+\int_{\Omega^e} K\Phi_{k}(\varphi)\,\qquad\varphi\in \mathcal{D}_0^{1,p}(\Omega^e)\,.
\end{equation}
Under our assumptions, $\mathcal{J}_k$ is well defined, continuous, coercive and strictly convex; indeed, $\Phi_k$ is convex and $|\Phi_k(s)|\leq|\Phi_k(0)|+k|s|$ for every $s\in\mathbb{R}$, so that, since $K\in L^1(\Omega^e)\cap L^{(p^*)'}(\Omega^e)$, the second integral in \eqref{Jk} is finite and continuous on $\mathcal{D}_0^{1,p}(\Omega^e)$, while $\varphi\mapsto\frac1p\int_{\Omega^e}H^p(\nabla\varphi)$ is strictly convex, continuous and, by \eqref{norm}, coercive.\\
Let $v\in W_{loc}^{1,p}\left(\Omega^e \right)$ be a supersolution to \eqref{problem} (see Definition~\ref{def:subsuper}). Notice that, by the definition of $\Phi_{k}$, it follows that $v$ is also a supersolution to the problem
\begin{equation}\label{supz}
-\Delta^{H}_{p} z=-K\Phi_{k}'(z).
\end{equation}

In what follows let us consider  the convex set
\begin{equation}\label{K}
	\mathcal K\,=\,\{\varphi\in \mathcal{D}_{0}^{1,p}\left(\Omega^e \right) \,:\, 0\leqslant \varphi\leqslant v\,\,\text{a.e. in}\,\,\Omega^e\}\,.
\end{equation}

Since $\mathcal K$ is a nonempty, closed and convex subset of $\mathcal{D}_0^{1,p}(\Omega^e)$, $\mathcal J_k$ admits a unique minimum point on $\mathcal K$. Denoting by $L^{\infty}_c(\Omega^e)$ the set of functions $g \in L^{\infty}(\Omega^e)$ such that $supp (g) \subset \subset \ \Omega^e$, let us prove the following preliminary lemmas. 
\begin{lemma}\label{ineq1} Let $w$ be a minimum of $\mathcal{J}_k$ on the set $\mathcal{K}$. If $\psi\in w+\left(\mathcal{D}_0^{1,p}(\Omega^e)\cap L^\infty_c(\Omega^e)\right)$ and $0\leqslant \psi\leqslant v$, then
	\begin{equation}
			\int_{\Omega^e }H ^{p-1} \left( \nabla w\right) \left\langle \nabla H \left(\nabla w \right), \nabla (\psi-w) \right\rangle  +\int_{\Omega^e}\,K \Phi_{k}'(w)(\psi-w)\geq 0.
	\end{equation}
\end{lemma}
\begin{proof}
For $t\in(0,1]$ we have $w+t(\psi-w)\in\mathcal K$, so that $\mathcal{J}_k(w)\leq\mathcal{J}_k(w+t(\psi-w))$. Since $\Phi_k\in C^1(\mathbb{R})$ with $|\Phi_k'|\leq k$ and $\psi-w\in L^\infty_c(\Omega^e)$, dividing by $t$ and letting $t\to0^+$ we get the thesis (see also \cite[Lemma~3.1]{EST}).
\end{proof}

Using the above lemma, we can prove the following

\begin{lemma}\label{ineq2}  Let $w$ be a minimum of $\mathcal{J}_k$ on the set $\mathcal{K}$. 
If $\psi \in C^\infty_c \left(\Omega^e\right)$ and $\psi\geq 0$ a.e. in $\Omega^e$, then
\begin{equation}
	\int_{\Omega^e }H ^{p-1} \left( \nabla w\right) \left\langle \nabla H \left(\nabla w \right), \nabla \psi \right\rangle+\int_{\Omega^e}\,K \Phi_{k}'(w)\psi \geq 0 
\end{equation}
\end{lemma}
\begin{proof}
	Let $\theta\in C^\infty_c(\mathbb{R})$ be with $\theta(t)\in[0,1]$ for $t\in\mathbb{R}$, $\theta (t)=1$
	for $t\in[-1,1]$ and $\theta (t)=0$ for $t\in(-\infty,-2]\cup[2,+\infty)$.	\\
	 Then, for any \emph{nonnegative} $\eta\in C^\infty_c(\Omega^e)$ we set
	\[
	\eta_h\,:=\,\theta\left(\frac{w}{h}\right) \,\eta,\qquad\quad \eta_{h,t}\,:=\,\min\{w+t\eta_h\,,\,v\}\,,
	\]
	with $h\geqslant 1$ and $t>0$. By the definition it follows that $\eta_{h,t}\leq v$; it is straightforward by the nonnegativity of $\eta$,  that $\eta_{h,t}\geq w$. Moreover $(\eta_{h,t}-w)\in \left(\mathcal{D}_0^{1,p}(\Omega^e)\cap L^\infty_c(\Omega^e)\right)$. By Lemma \ref{ineq1} it follows that
	
	\begin{equation}\label{eq w}
			\int_{\Omega^e }H ^{p-1} \left( \nabla w\right) \left\langle \nabla H \left(\nabla w \right), \nabla (\eta_{h,t}-w) \right\rangle+ \int_{\Omega^e}\,K \Phi_{k}'(w)(\eta_{h,t}-w)\geq 0.
	\end{equation}
	Consequently,  using Lemma \ref{estimates} and $\eqref{eq w}$ we obtain
	\begin{eqnarray*}
			&C_{1}&\int_{\Omega^e }\left( \left\vert \nabla
			\eta_{h,t}\right\vert +\left\vert \nabla w\right\vert \right) ^{p-2}\left\vert
			\nabla \left( \eta_{h,t}-w\right) \right\vert^{2}+\int_{\Omega^e}\,K \left(\Phi_{k}'(\eta_{h,t})-\Phi_{k}'(w)\right)(\eta_{h,t}-w) \\
		&\leq& \int_{\Omega^e } \left\langle H ^{p-1} \left( \nabla \eta_{h,t}\right)\nabla H \left(\nabla \eta_{h,t} \right)- H ^{p-1} \left( \nabla w\right)\nabla H \left(\nabla w \right), \nabla \left( \eta_{h,t}-w\right) \right\rangle \\
			&&+ \int_{\Omega^e}\,K \left(\Phi_{k}'(\eta_{h,t})-\Phi_{k}'(w)\right)(\eta_{h,t}-w)\\
			&=& 	\int_{\Omega^e }H ^{p-1} \left( \nabla \eta_{h,t}\right) \left\langle \nabla H \left(\nabla \eta_{h,t} \right), \nabla (	\eta_{h,t}-w) \right\rangle+\int_{\Omega^e}\,K\Phi_{k}'(\eta_{h,t})(\eta_{h,t}-w) \\
			&&- \int_{\Omega^e}H ^{p-1} \left( \nabla w\right) \left\langle \nabla H \left(\nabla w \right), \nabla (\eta_{h,t}-w) \right\rangle-\int_{\Omega^e}\,K\Phi_{k}'(w)(\eta_{h,t}-w)\\
			&\leq& \int_{\Omega^e }H ^{p-1} \left( \nabla \eta_{h,t}\right) \left\langle \nabla H \left(\nabla \eta_{h,t} \right), \nabla (	\eta_{h,t}-w) \right\rangle+\int_{\Omega^e}\,K\Phi_{k}'(\eta_{h,t})(\eta_{h,t}-w)\\
		&=&  \int_{\Omega^e }H ^{p-1} \left( \nabla \eta_{h,t}\right) \left\langle \nabla H \left(\nabla \eta_{h,t} \right), \nabla (	\eta_{h,t}-w-t\eta_{h}) \right\rangle+\int_{\Omega^e}\,K\Phi_{k}'(\eta_{h,t})(\eta_{h,t}-w-t\eta_{h})\\
			&&+  t\int_{\Omega^e }H ^{p-1} \left( \nabla \eta_{h,t}\right) \left\langle \nabla H \left(\nabla \eta_{h,t} \right), \nabla \eta_{h} \right\rangle+t\int_{\Omega^e}\,K\Phi_{k}'(\eta_{h,t})\eta_{h}.
\end{eqnarray*}
Since 
\begin{equation*}
	\nabla (\eta_{h,t}-w-t\eta_{h})=\,
	\begin{cases}
		\nabla (v-w-t\eta_{h})\qquad\text{if}\quad w+t\eta_h\geq v\\
0\qquad \quad \qquad otherwise
	\end{cases}
\end{equation*}
we have
\begin{eqnarray*}
&C_{1}&\int_{\Omega^e }\left( \left\vert \nabla
			\eta_{h,t}\right\vert +\left\vert \nabla w\right\vert \right) ^{p-2}\left\vert
			\nabla \left( \eta_{h,t}-w\right) \right\vert^{2}+\int_{\Omega^e}\,K \left(\Phi_{k}'(\eta_{h,t})-\Phi_{k}'(w)\right)(\eta_{h,t}-w) \\
&\leq&\int_{\Omega^e }H ^{p-1} \left( \nabla \eta_{h,t}\right) \left\langle \nabla H \left(\nabla \eta_{h,t} \right), \nabla (	\eta_{h,t}-w-t\eta_{h}) \right\rangle+\int_{\Omega^e}\,K\Phi_{k}'(\eta_{h,t})(\eta_{h,t}-w-t\eta_{h})\\
			&&+  t\int_{\Omega^e }H ^{p-1} \left( \nabla \eta_{h,t}\right) \left\langle \nabla H \left(\nabla \eta_{h,t} \right), \nabla \eta_{h} \right\rangle+t\int_{\Omega^e}\,K\Phi_{k}'(\eta_{h,t})\eta_{h}\\
			&=& \int_{\Omega^e }H ^{p-1} \left( \nabla v\right) \left\langle \nabla H \left(\nabla v\right), \nabla (	\eta_{h,t}-w-t\eta_{h}) \right\rangle+\int_{\Omega^e}\,K\Phi_{k}'(v)(\eta_{h,t}-w-t\eta_{h})\\
	&&+  t\int_{\Omega^e }H ^{p-1} \left( \nabla \eta_{h,t}\right) \left\langle \nabla H \left(\nabla \eta_{h,t} \right), \nabla \eta_{h} \right\rangle+t\int_{\Omega^e}\,K\Phi_{k}'(\eta_{h,t})\eta_{h}.
	\end{eqnarray*}
Using the fact that $v$ is supersolution to \eqref{supz} and observing also that $\eta_{h,t}-w- t\eta_h\leq 0$, by a  density argument we deduce
$$
 \int_{\Omega^e }H ^{p-1} \left( \nabla v\right) \left\langle \nabla H \left(\nabla v\right), \nabla (	\eta_{h,t}-w-t\eta_{h}) \right\rangle+\int_{\Omega^e}\,K\Phi_{k}'(v)(\eta_{h,t}-w-t\eta_{h})\leq 0
$$
hence
	\begin{eqnarray*}
		&C_{1}&\int_{\Omega^e }\left( \left\vert \nabla
		\eta_{h,t}\right\vert +\left\vert \nabla w\right\vert \right) ^{p-2}\left\vert
		\nabla \left( \eta_{h,t}-w\right) \right\vert^{2}+\int_{\Omega^e}\,K \left(\Phi_{k}'(\eta_{h,t})-\Phi_{k}'(w)\right)(\eta_{h,t}-w) \\
		&\leq&  t\int_{\Omega^e }H ^{p-1} \left( \nabla \eta_{h,t}\right) \left\langle \nabla H \left(\nabla \eta_{h,t} \right), \nabla \eta_{h} \right\rangle+t\int_{\Omega^e}\,K\Phi_{k}'(\eta_{h,t})\eta_{h}.
	\end{eqnarray*}
Since $\Phi_k'$ is nondecreasing and $\eta_{h,t}\geq w$, the left-hand side of the last inequality is nonnegative; hence, dividing by $t>0$, we get
$$
\int_{\Omega^e }H ^{p-1} \left( \nabla \eta_{h,t}\right) \left\langle \nabla H \left(\nabla \eta_{h,t} \right), \nabla \eta_{h} \right\rangle+\int_{\Omega^e}\,K\Phi_{k}'(\eta_{h,t})\eta_{h}\geq 0.
$$
	Finally, passing to the limit for $t\rightarrow 0^{+}$, exploiting the Lebesgue dominated convergence theorem, we deduce that
	\begin{eqnarray*}
			\int_{\Omega^e }H ^{p-1} \left( \nabla w\right) \left\langle \nabla H \left(\nabla w \right), \nabla \eta_{h} \right\rangle+\int_{\Omega^e}K\Phi_{k}'(w)\eta_{h}\geq 0.
	\end{eqnarray*}
	The claim follows passing to the limit for $h\rightarrow +\infty$ in the last inequality.
\end{proof}

\begin{theorem}\label{comparison}
	Let $\gamma>0$ and let $K$ be a positive function satisfying \eqref{hyp:K}. If $u,v\in W^{1,p}_{loc}(\Omega^e)$ are a subsolution and a supersolution to \eqref{problem}, respectively, then $u\leq v$ a.e. in $\Omega^e$.
\end{theorem}
\begin{proof}
	Fix $\varepsilon>0$ and $k>\varepsilon^{-\gamma}$, and let us consider again $w\in \mathcal{D}_0^{1,p}(\Omega^e)$, $w\geqslant 0$ a.e. in $\Omega^e$, the minimum of the functional $\mathcal{J}_{k}$ on the set $\mathcal{K}$, defined in \eqref{Jk} and \eqref{K}, respectively. 

Since $(u-\varepsilon)^+\in\mathcal{D}_0^{1,p}(\Omega^e)$ (see Definition~\ref{def:subsuper}) and $w\geq0$, we deduce that $(u-\varepsilon-w)^+=\big((u-\varepsilon)^+-w\big)^+\in\mathcal{D}_0^{1,p}(\Omega^e)$.	Defining  
	$$T_\tau (s):=\min\{s,\tau\},\mbox{ if } s\geq 0\quad \mbox{and} \quad T_\tau (-s)=-T_\tau (s) \mbox{ if } s< 0,$$
	by a density argument and Lemma \ref{ineq2} it follows
	\begin{equation}\label{supermin}
		\int_{\Omega^e }H ^{p-1} \left( \nabla w\right) \left\langle \nabla H \left(\nabla w \right), \nabla T_\tau\left((u-\varepsilon-w)^+\right) \right\rangle+\int_{\Omega^e}\,K \Phi_{k}'(w) T_\tau\left((u-\varepsilon-w)^+\right)\geq 0.
	\end{equation}
\noindent Let us consider $\varphi_n\in C^\infty_c(\Omega^e)$ such that $\varphi_n\rightarrow (u-\varepsilon-w)^+$ in $\mathcal{D}_0^{1,p}(\Omega^e)$ and set 
	\begin{equation*}
		\tilde\varphi_{\tau,n}:=T_\tau (\min\{(u-w-\varepsilon)^+,\varphi_n^+\}). 
	\end{equation*}
in such a way that $\tilde\varphi_{\tau,n}\in \mathcal{D}_0^{1,p}(\Omega^e)\cap L^\infty_c(\Omega^e)$. Hence, since $u$ is a subsolution to \eqref{problem}, taking $\tilde\varphi_{\tau,n}$ as test function, we get
	\begin{equation*}
	\int_{\Omega^e}H ^{p-1} \left( \nabla u\right) \left\langle \nabla H \left(\nabla u \right), \nabla \tilde\varphi_{\tau,n} \right\rangle \leq
	\int_{\Omega^e}\,\frac{K}{u^\gamma}\,\tilde\varphi_{\tau,n}.
	\end{equation*}
	Passing to the limit as $n$ tends to infinity in the above inequality (note that $u>\varepsilon$ on the support of $\tilde\varphi_{\tau,n}$, so that $\nabla u=\nabla(u-\varepsilon)^+\in L^p(\Omega^e)$ there, and $0\leq Ku^{-\gamma}\tilde\varphi_{\tau,n}\leq\varepsilon^{-\gamma}K(u-w-\varepsilon)^+\in L^1(\Omega^e)$), it is easy to deduce that
	\begin{equation}\label{sub}
		\int_{\Omega^e}H ^{p-1} \left( \nabla u\right) \left\langle \nabla H \left(\nabla u \right), \nabla T_\tau\left((u-w-\varepsilon)^+\right) \right\rangle \leq \int_{\Omega^e}\frac{K}{u^\gamma}T_\tau\left((u-w-\varepsilon)^+\right).
	\end{equation}
Subtracting \eqref{supermin} from \eqref{sub} and using Lemma \ref{estimates} we have
\begin{eqnarray*}
			&C_{1}& \int_{\Omega^e} \left( \left\vert \nabla
			u\right\vert +\left\vert \nabla w\right\vert \right) ^{p-2}\left\vert\nabla T_\tau\left((u-w-\varepsilon)^+\right)\right\vert^2 \\
			&\leq& \int_{\Omega^e}\left\langle H ^{p-1} \left( \nabla u\right)\nabla H \left(\nabla u\right)- H ^{p-1} \left( \nabla w\right)\nabla H \left(\nabla w \right), \nabla T_\tau\left((u-w-\varepsilon)^+ \right) \right\rangle \\
			&\leq& \int_{\Omega^e} K \left( \frac{1}{u^{\gamma}}+\Phi_{k}'(w)\right) T_\tau \left( (u-w-\varepsilon)^+  \right).
	\end{eqnarray*}

	Now, since $k>\varepsilon^{-\gamma}$ and $u>\varepsilon$ on the support of $T_\tau((u-w-\varepsilon)^+)$, there we have that $\Phi_{k}'(u)=-u^{-\gamma}$; hence, since $\Phi_k'$ is nondecreasing, we obtain
	\begin{eqnarray*}
			&C_{1}& \int_{\Omega^e} \left( \left\vert \nabla
			u\right\vert +\left\vert \nabla w\right\vert \right) ^{p-2}\left\vert\nabla T_\tau\left((u-w-\varepsilon)^+\right)\right\vert^2 \\
			&\leq& \int_{\Omega^e}\left\langle H ^{p-1} \left( \nabla u\right)\nabla H \left(\nabla u\right)- H ^{p-1} \left( \nabla w\right)\nabla H \left(\nabla w \right), \nabla T_\tau\left((u-w-\varepsilon)^+ \right) \right\rangle \\
			&\leq& \int_{\Omega^e} K\left(-\Phi_{k}'(u)+\Phi_{k}'(w)\right)
			T_\tau\left((u-w-\varepsilon)^+\right)\leq 0.
	\end{eqnarray*}
	By the arbitrariness  of $\tau>0$,  we deduce that $u\leq w+\varepsilon \leq v+\varepsilon$ a.e. in $\Omega^e$ and the thesis follows letting $\varepsilon\rightarrow 0$.	
\end{proof}

As a direct consequence of Theorem \ref{comparison} we obtain the following result.
\begin{theorem}\label{uniqueness 1}
	Let $\gamma>0$ and let $K$ be a positive function satisfying \eqref{hyp:K}. Then the solution $u$, in the sense of Definition \rev{\ref{def:sol}}, to problem \eqref{problem} is unique.
\end{theorem}

\begin{proof}
It suffices to observe that every solution of \eqref{problem} in the sense of Definition~\ref{def:sol} is both a subsolution and a supersolution in the sense of Definition~\ref{def:subsuper}, and to apply Theorem~\ref{comparison} twice.
\end{proof}

\begin{proof}[Proof of Theorem~\ref{thm:uniqueness}]
It follows from Theorems~\ref{comparison} and \ref{uniqueness 1}.
\end{proof}

\end{document}